\documentclass[12pt]{article}

\usepackage{enumerate}
\usepackage{graphicx}
\usepackage{amsmath}
\usepackage{amssymb}
\usepackage{epsfig}
\usepackage{amsthm}
\usepackage{amscd}
\usepackage[latin1]{inputenc}
\def\r{\mathbb{R}}
\def\b{\mathbb{B}}

\def\n{\mathbb{N}}
\def\c{\mathbb{C}}
\def\s{\mathbb{S}}

\newcommand{\df}{ \stackrel{\rm def}{=}}

\newcommand{\rth}{\r^3}

\newcommand{\dD}{d_{\partial {\cal D}}}

\newtheorem{theorem}{Theorem}

\begin{document}

\title{\sc Bounded domains which are universal for minimal surfaces}
\author{\\ Francisco Martín\thanks{Research partially supported by MEC-FEDER Grant no. MTM2004 - 00160.
} \and  \\ William H. Meeks, III \thanks{This material is based upon
 work for the NSF under Award No. DMS - 0405836. Any opinions, findings, and conclusions or recommendations
 expressed in this publication are those of the authors and do not
 necessarily reflect the views of the NSF.}
 \and \\ Nicolai Nadirashvili}
\date{\today}
\maketitle

\begin{abstract}  We construct open domains in $\rth$ which do not
admit complete properly immersed minimal surfaces with an annular
end.  These domains can not be smooth by a recent result of Martín
and  Morales \cite{marmor2}.

\vspace{.2cm} \noindent   {\it 2000 Mathematics Subject
Classification.} Primary 53A10; Secondary 49Q05, 49Q10, 53C42.

\end{abstract}

\section{Introduction}
 The main goal of this paper is to construct bounded open domains in
 $\r^3$ which do not contain any complete properly immersed minimal
 surfaces with at least one annular end. It is our belief that these
  open domains are in fact universal according to the following
   definition: A connected region of space which is open or the
    closure of an open set is {\it universal for minimal surfaces},
     if every complete properly immersed minimal surface in the region
     is recurrent for Brownian motions. In particular, a bounded domain
     is universal if and only if it contains no complete properly immersed
     minimal surfaces.
\begin{theorem} \label{th:first}
Let ${\cal D}$ be any bounded open domain in $\r^3$. Then there exists a
proper countable collection $\cal F$ of pairwise disjoint horizontal
simple closed curves in $\cal D$ such that the complementary domain
 $\widetilde{\cal D}={\cal D} - {\cal F}$ is universal for minimal surfaces
  with at least one annular end. In particular, any complete immersed minimal
   surface of finite genus in $\widetilde{\cal D}$ must have an uncountable number of ends.
\end{theorem}
 The construction of the domains $\widetilde{\cal D}$ that appear in the
 above theorem are motivated by a related unpublished example of the third
  author. We will explain a variant of his original example at the end of Section 2.

Interest in results like Theorem \ref{th:first} dates back to an earlier
 question by Calabi. Calabi asked whether or not it is possible for a complete
  minimal surface in $\r^3$ to be contained in the ball $B=\{x \in \r^3 \; | \; \|x\|<1\}.$
In \cite{na1}, Nadirashvili constructed a complete minimal surface in $B$. After
 Nadirashvili negative solution to Calabi's question, Martín and
 Morales \cite{marmor1} proved that there exist complete properly
 immersed minimal disks in $B$. Recently \cite{marmor3}, they improved
 on their original techniques and were able to show that every bounded
 domain with $C^{2,\alpha}$-boundary admits a complete properly immersed minimal disk
 whose boundary limit set is close to a prescribed simple closed curve on the
  boundary of the domain. In contrast to these existence results for complete
  properly immersed minimal disks in bounded domains, Colding and Minicozzi \cite{cm35}
   recently proved that any complete embedded minimal surface in $\rth$ is properly
    embedded in $\rth$. By results of Meeks and Rosenberg, \cite{mr6, mr8}, any
     properly embedded minimal surface of finite topology in $\rth$ is recurrent
     for Brownian motion. Hence, every domain in $\rth$ is universal for embedded
      minimal surfaces of finite topology. Finally, we remark that
       Collin, Kusner, Meeks and Rosenberg \cite{ckmr1} proved that any properly
        immersed minimal surface with boundary in a closed convex domain in $\rth$
         has full harmonic measure on its boundary.

At the end of Section 2, we give an estimate for the growth of the absolute curvature
function $|K_M|$  for any complete properly immersed minimal surface $M$ in a smooth
 bounded domain ${\cal D} \subset \rth$ in terms of the distance function
 $d_{\partial {\cal D}}$ of $M$ to $\partial {\cal D}$. This estimate implies
 the function $|K_M| \,d^2_{\partial {\cal D}}$ is never bounded.

\section{Proof of Theorem 1}

Let $\cal D$ be an open connected bounded set of $\rth$ and let $\overline{\cal D}$.
Without loss of generality we may assume that $\overline{\cal D}$ is contained in the closed slab
$$\left\{(x_1,x_2,x_3) \in \r^3 \; | \; 0 \leq x_3 \leq 1 \right\}$$
and $\overline{\cal D}$ contains points at heights $0$ and $1$.

For $t \in (0,1)$, let $P_t$ denote the horizontal plane at height $t$. Let $C_t={\cal D} \cap P_t$,
which consists of a collection $\{C_{t,i}\}_{i \in I_t}$ of connected components, for some countable
indexing set $I_t.$ For each $t$ and for each $i \in I_t$, choose an exhaustion of $C_{t,i}$ by
smooth compact domains $C_{t,i,k}$, $k \in \n$, and where
$C_{t,i,k} \subset C_{t,i,k+1}, \quad \forall k \in \n.$
Finally, let $\displaystyle C_{t}(k) \df \bigcup_{i \in I_t} C_{t,i,k}.$
\begin{figure}[htb]
    \begin{center}
        \includegraphics[width=0.94\textwidth]{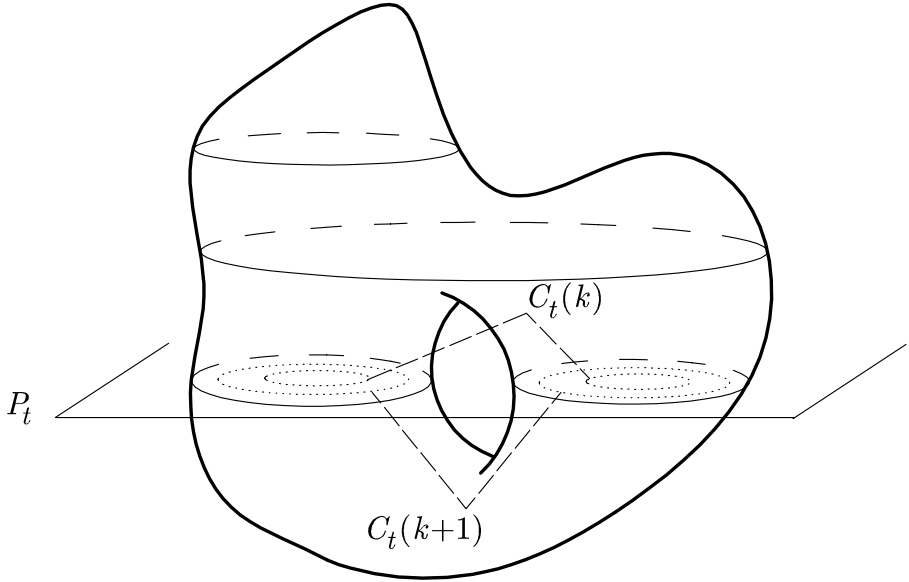}
    \end{center}
    \caption{\small The domain $\cal D$ and the sets $C_t(k).$}
    \label{fig:lospe}
\end{figure}

Now consider the following sequence of ordered rational numbers:
$$Q=\left\{\frac 12 \,,
 \,\frac 13\,,\, \frac 23\,,\, \frac 14\,,\,
  \frac 24\,,\, \frac 34\,,\, \ldots\,,\, \frac 1n\,,\,
   \frac 2n\,, \,\ldots\,,\, \frac{n-1}n\,,\, \ldots \right\}.$$
Let $t_k$ the $k$-th rational number in $Q$. Define $\cal F$ to be the collection
of boundary curves to all of the domains $\displaystyle \bigcup_{k \in \n} C_{t_k}(k),$
and define $\widetilde{\cal D} \df {\cal D} - {\cal F}.$

Suppose that $f:M \rightarrow \widetilde{\cal D}$ is a complete properly immersed minimal
surface with an annular end $E$ and we will obtain a contradiction. First, note that the
limit set $L(E)$ of $E$ is a connected set contained in $\partial \widetilde{\cal D},$ where
$$\partial \widetilde{\cal D}=\overline{\,\widetilde{\cal D}\,}-\widetilde{\cal D}=\partial{\cal D}
 \cup {\cal F}.$$

Our initial goal is to prove that $x_ 3|_{L(E)}$ is constant, from which we will easily obtain a contradiction.

If $L(E)$ intersects one of the horizontal curves $C$ in $\cal F$, then $L(E) \subset C$
(recall that $L(E)$ is connected) and we have proved that  $x_ 3|_{L(E)}$ is constant. So,
suppose that $p \in L(E) \subset \partial {\cal D}.$ If $x_ 3|_{L(E)}$ is not constant, then
 there exists a point $q\in L(E)$ with $x_3(p) \neq x_3(q)$. Choose a positive rational number
 $t$ which lies between $x_3(p)$ and $x_3(q)$. Notice that $t$ can be represented by an infinite
 subsequence $\{t_{k_1} \, , \, t_{k_2} \, , \,\ldots \, , \,  t_{k_n} \, , \, \ldots \} \subset Q$.
  Since the plane $P_t$ separates $p$ and $q$,  for every subend $E' \subset E$, $P_t \cap E'$
  is nonempty. On the other hand, the subdomains $C_t(k_n)$ give a compact exhaustion
  to $P_t\cap {\cal D}$ with boundaries disjoint from $E$. Therefore, every component of
   $P_t \cap E$ is compact. Since  $P_t \cap E$ is noncompact, then there exist a pair of
    disjoint simple closed curves in $P_t \cap E \subset E$ which bound a compact domain
    in $E$, since $E$ is an annulus. But then the harmonic function $x_3$ restricted to this
     domain has an interior maximum or minimum which is impossible. This contradiction proves
     that $x_ 3|_{L(E)}$ {\it is constant.} Let $a$ denote this constant.
\begin{figure}[htbp]
    \begin{center}
        \includegraphics[width=0.75\textwidth]{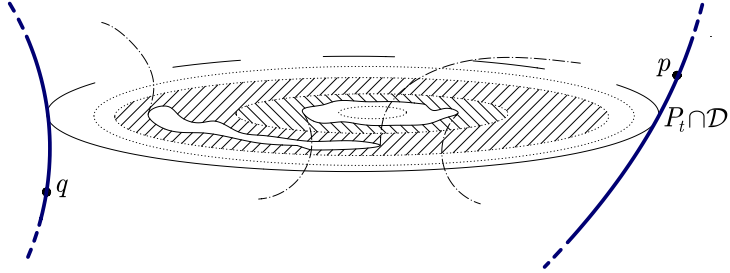}
    \end{center}
    \caption{\small The subdomains $C_t(k_n)$ give a compact exhaustion to $P_t\cap {\cal D}$
     with boundaries disjoint from $E$. Therefore, every component of $P_t \cap E$ is compact.
     Since  $P_t \cap E$ is noncompact, then there exist a pair of disjoint simple closed curves
      in $P_t \cap E \subset E$ which bound a compact domain in $E$}
    \label{fig:Drawing-3}
\end{figure}

Our next step consists of proving that if $x_ 3|_{L(E)}$  is
constant, then the minimal  immersion $f:M \rightarrow
\widetilde{\cal D}$ is incomplete, which is contrary to our
assumptions. Indeed, consider a conformal parameterization of the
end $E$ by the annulus $A=\{z \in \c \; | \; r \leq |z|<1\} \subset
\c$, for some $0<r<1$. Since $x_3$ is a
 {\it bounded } harmonic function defined on  $A$, then by Fatou's theorem $x_3$ has
 radial limit a.e. in $\s^1=\{z \in \c \;|\; |z|=1\}.$ Furthermore, the function $x_3$
 is determined by the Poisson integral of its radial limits (see for instance \cite{conway2}.)
 Since the limit $\lim_{\rho\to 1}x_3(\rho \theta)=a$, at almost every point $\theta$ in $ \s^1$,
 then $x_3$ admits a regular extension to $\overline{A}$. In particular, $\| \nabla x_3 \|$ is bounded
  in $A$. On the other hand, as $x_2$ is also a bounded harmonic function, then a result by
   Bourgain \cite[Theorem 2]{bourgain1} asserts that the set
$${\cal S}=\left\{ \theta \in \s^1 \; | \;  \int_r^1 \|\nabla x_2(\rho \,\theta)\|\,d\rho<+\infty \right\} $$
has Hausdorff dimension 1, in particular $\cal S$ is nonempty. Moreover, for a conformal minimal
immersion it is well known \cite{os1} that $\|\nabla x_1 \| \leq \| \nabla x_2 \|+\| \nabla x_3 \|.$

Hence, as a consequence of all these facts,  if $\theta$ is a point in $\cal S$ then
 $$\int_r^1 \sqrt{\| \nabla x_1(\rho \, \theta)\|^2+\| \nabla x_2(\rho \,
  \theta)\|^2+\| \nabla x_3(\rho \, \theta)\|^2} \;d \rho < \infty,$$ which means that
   the divergent curve $f(\rho \, \theta)$, $\rho \in (r,1)$, has finite length, and so
    $f$ is not complete. This contradiction proves the theorem.
\vskip .5cm

We now explain a modification of the original unpublished example of
Nadirashvili which motivated our construction of the domains
$\widetilde{\cal D}$ given in Theorem \ref{th:first}.

Let $\cal D$ be the open cube:
$${\cal D}=(-1,1) \times (-1,1) \times (-1,1).$$
Let $F_1=\{1\} \times [-1,1] \times [-1,1]$,  $F_2=[-1,1] \times \{1\} \times [-1,1]$ and
 $F_3= [-1,1]\times [-1,1] \times \{1\}$ be the three coordinate faces of $\cal D$. Let
  $S_i \df \partial F_i$ be the related boundary square curves. As in the construction of
   the domains in Theorem \ref{th:first}, we need to define a countable proper collection
    $\cal F$ of planar simple closed curves in the cube $\cal D$, so that ${\cal D}- {\cal F}$
    admits no complete properly immersed minimal surfaces with an annular end.
\begin{figure}[htbp]
    \begin{center}
        \includegraphics[width=0.75\textwidth]{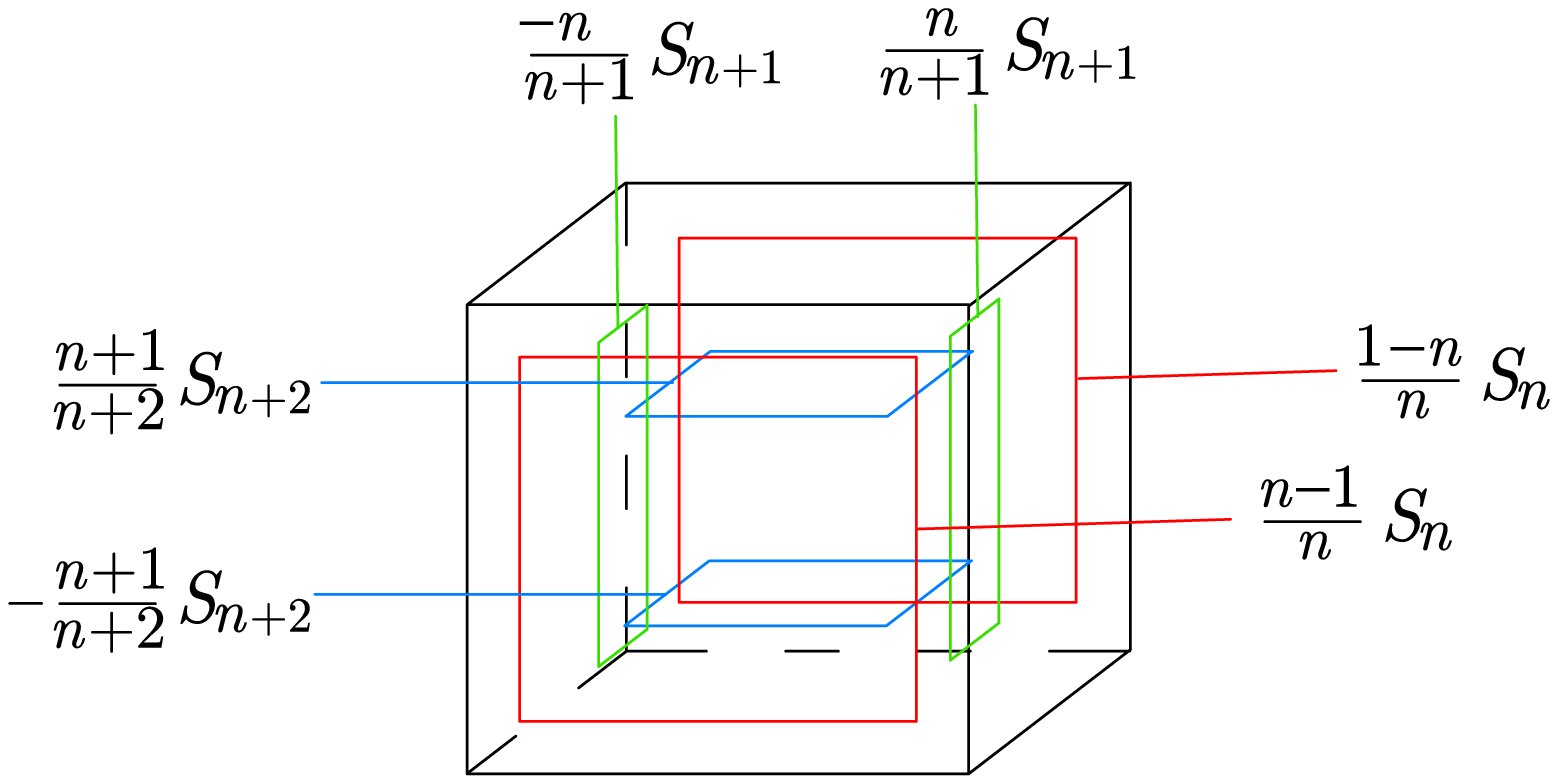}
    \end{center}
    \caption{The cube $\cal D$.}
    \label{fig:Drawing-4}
\end{figure}

For a real number $\lambda$, let $\lambda \, S_i=\{ \lambda \,(x_1,x_2,x_3) \; |
\; (x_1,x_2,x_3) \in S_i\},$ for $i=1,2,3.$ Let $\cal F$ be the collection of curves
$$\left\{ \tfrac{n-1}{n} \, S_{n ({\rm mod} 3)}\, , \; \tfrac{1-n}{n} \, S_{n ({\rm mod} 3)}
\right\}_{n \in \n}.$$
Then, a small modification of the arguments given in the proof of Theorem \ref{th:first}
implies that one of the coordinate functions $(x_1,x_2,x_3)$ restricted to the limit set
of an annular end of a complete immersed minimal surface in $\widetilde{\cal D} \df {\cal D}-{\cal F}$
is constant. As in the proof of Theorem \ref{th:first}, we obtain a contradiction.
\vskip .5cm

Finally, we explain that if $M$ is a complete properly immersed minimal surface in a convex or
 smooth bounded domain, then the function $|K_M| \, {\dD}^2$ is not bounded.

We proceed by contradiction. Assume there exists a constant $C>0$ so that
$|K_M| \, d_{\partial D}\leq C.$
Since the convex hull of a nonflat complete minimal surface with bounded curvature in
$\rth$ is all of $\rth$ \cite{xa5}, the curvature function is not bounded in $M$. Thus,
 take an arbitrary sequence of points $q_n \in M$ such that $|K_M(q_n)| \geq n^2$. Let
  $p_n' \in M \cap B_M(q_n,1)$ be a maximum of $$h_n=|K_M| \,d_M(\cdot, \partial B_M(q_n,1))^2,$$
where $d_M$ denotes the intrinsic distance of $M$ and $B_M(q_n,1)$ means the intrinsic
 ball centered at $q_n$ with radius 1.

We label $\lambda'_n=\sqrt{|K_M(p_n')|}.$ Notice that:
\begin{multline*} \lambda'_n \geq \lambda_n' \, d_M(p_n', \partial B_M(q_n,1))=
\sqrt{h_n(p_n')} \geq \\
\sqrt{h_n(q_n)}=\sqrt{|K_M(q_n)|}=n.\end{multline*} Fix $t>0$.
Notice that the sequence of extrinsic balls $$\left\{ \lambda_n' \b
\left(p_n' , \frac{t}{\lambda_n'} \right)\right\}_{n \in \n}$$
converges to the ball $\b(t)$, where we
 have indentified $p'_n$ with $\vec 0.$ Similarly, we can consider  $\{\lambda_n' B_M(p_n',t/\lambda_n')
 \}_{n \in \n}$ as a sequence of minimal surfaces with boundary, passing through $\vec 0$ with curvature
  $-1$ at the origin. From our assumption, we know that $D_n(t)=\partial {\cal D} \cap \b
  \left(p_n' ,\frac{t}{\lambda_n'} \right)$ is nonempty, for any $t>C.$

We assert that the curvature of these minimal surfaces with boundary is uniformly bounded.
 Indeed, pick a point $z$ in $B_M(p_n',t/\lambda_n').$ Then we have
\begin{equation} \label{eq:ca}
\frac{\sqrt{|K_M(z)|}}{\lambda'_n}=\frac{\sqrt{h_n(z)}}{\lambda_n' \, d_M(z,\partial B_M(q_n,1))}
\leq \frac{d_M(p_n',\partial B_M(q_n,1))}{d_M(z,\partial B_M(q_n,1))}
\end{equation}
By the triangle inequality, one has $$d_M(p_n',\partial B_M(q_n,1))\leq \frac{t}{\lambda_n'}+
d_M(z,\partial B_M(q_n,1)),$$
and so
\begin{multline*}\frac{d_M(p_n',\partial B_M(q_n,1))}{d_M(z,\partial B_M(q_n,1))} \leq 1+
\frac{t}{\lambda_n' \, d_M(z,\partial B_M(q_n,1))} \leq \\
1+\frac{t}{\lambda'_n \, \left( d_M(p_n',\partial B_M(q_n,1))-\frac{t}{\lambda_n'}\right)}
 \leq 1+ \frac{t}{n-t},\end{multline*}
which tends to $1$ as $n \to \infty.$

After extracting a subsequence, it follows that $\lambda_n'\, B_M(p_n',\frac{t}{\lambda_n'})$
converge smoothly to a minimal surface $M_\infty(t)$ contained in $\b(t).$ Since
$\displaystyle \lim_{n \to \infty }\lambda_n' =+\infty$, then $\lambda_n' \, D_n(t)$
converges either to a plane in the case that $\cal D$ is a regular domain or to the boundary of
 a convex body if $\cal D$ is a convex domain. In any case, $M_\infty(t)$ is contained in one of
  the halfspaces determined by the plane, or in the interior of the convex body. Note that
   $M_\infty=\cup_{t\geq C} M_\infty(t)$ is a complete nonflat minimal surface. By construction,
   $M_\infty$ has bounded curvature and is contained in a convex domain which is not $\r^3$. But
    this is contrary to the aforementioned result by Xavier. This contradiction proves that
    $|K_M| \, {\dD}^2$ is not bounded.

\center{William H. Meeks, III at bill@math.umass.edu}\\
Mathematics Department at the University of Massachusetts\\
Amherst, MA 01003, USA.
 \center{Francisco Martín at fmartin@ugr.es }\\
Departamento de Geometría y Topología. Universidad de Granada\\
18071 Granada, Spain. 
   \center{ Nicolai Nadirashvili at nicolas@cmi.univ-mrs.fr \\
   CNRS, LATP, CMI, 39, rue Joliot-Curie \\
   13453 Marseille cedex 13, France}

\end{document}